\documentclass[leqno]{amsart}
\usepackage{amssymb}

\newtheorem{thm}{Theorem}[section]
\newtheorem{prop}[thm]{Proposition}
\newtheorem{lemma}[thm]{Lemma}

\theoremstyle{remark}
\newtheorem{ex}[thm]{Example}
\newtheorem{rmk}[thm]{Remark}
\theoremstyle{definition}
\newtheorem{defn}[thm]{Definition}

\numberwithin{equation}{section}

\renewcommand{\L}{\mathcal{L}}
\newcommand{\X}{\mathcal{X}}
\renewcommand{\H}{\mathcal{H}}

\newcommand{\vf}{\varphi}

\newcommand{\ox}{\otimes}
\newcommand{\wt}{\widetilde}
\renewcommand{\:}{\colon}
\newcommand{\IP}{{\mathbb P}}
\renewcommand{\div}{\mathrm{div}}
\newcommand{\Pic}{\mathrm{Pic}}
\renewcommand{\Im}{\mathrm{Im}}
\newcommand{\Ker}{\mathrm{Ker}}
\newcommand{\lra}{\longrightarrow}
\renewcommand{\O}{\mathcal{O}}
\def\Hilb{\operatorname{Hilb}}
\def\low{\operatorname{low}}
\def\high{\operatorname{high}}

\begin{document}

\title{Abel maps and limit linear series}
\author{Eduardo Esteves}
\author{Brian Osserman}
\begin{abstract} We explore the relationship between limit linear series
and fibers of Abel maps in the case of curves with two smooth components
glued at a single node. 
To
an $r$-dimensional 
limit 
linear series satisfying 
a certain exactness property (weaker than the refinedness property of 
Eisenbud and Harris) we associate a closed subscheme of the appropriate 
fiber of the Abel map. We then describe this closed subscheme explicitly, 
computing its Hilbert polynomial and showing that it is Cohen--Macaulay 
of pure dimension $r$. We show that this construction is also compatible
with one-parameter smoothings.
\end{abstract}

\thanks{
The first author was supported by 
CNPq, Proc.~303797/2007-0 and 473032/2008-2, and FAPERJ, 
Proc.~E-26/102.769/2008 and E-26/110.556/2010. This project 
was initiated while the authors were visiting MSRI for the 2009 program 
in Algebraic Geometry and concluded when the second author was visiting 
IMPA. The authors thank both institutes for the support given.
}
\maketitle

\section{Introduction}

The classical theory of linear series on smooth curves is closely
related to 
that of 
Abel maps and their fibers, which consist precisely
of complete linear series. This relationship also amplifies the relationship
between linear series and (families of) effective divisors. For 
(singular)
curves of compact type, Eisenbud and Harris \cite{e-h1}
developed the theory of limit linear series as
an analogue of linear series, while Coelho and Pacini \cite{c-p1}
have 
studied
Abel maps. However, the relationship
between these two concepts is far murkier than in the smooth case.
On the side of limit linear series, there is no obvious concept of a 
complete limit linear series, nor of families of divisors associated to
a limit linear series. On the other hand, fibers of 
Abel maps are not
very well behaved: for instance, they are in general not even 
equidimensional. 

Our aim is to relate limit linear series to fibers of Abel
maps via the 
definition
of limit linear series 
and construction of their moduli space 
in \cite{os8}.
For the sake of simplicity, we restrict our attention
to the case treated in {\it loc.~cit.},
which is that of a curve
$X$ with two smooth components glued together at a single node.

There is an open subset of the moduli space of limit linear series
on $X$ consisting of ``exact'' limit linear series (see Definition
\ref{defn:exact} below). These contain in particular all limits of linear
series on the generic fiber in a regular smoothing family; see Section
\ref{sec:limits}. If
$\mathfrak g$ is an exact limit linear series of dimension $r$ with
underlying line bundle $\L$, we
construct a closed subscheme $\IP(\mathfrak g)$ of $A_d^{-1}(\L)$, the
corresponding fiber of the $d$th Abel map. This subscheme is by definition
reduced, and we show in Theorem \ref{main} that $\IP(\mathfrak g)$ is 
connected and Cohen--Macaulay, of dimension $r$, with the same Hilbert 
polynomial as $\IP^r$. We also show in Theorem \ref{thm:flat-limit} that 
if $\mathfrak g$ is the limit of a $\mathfrak g^r_d$ on the generic
fiber of a one-parameter regular
smoothing of $X$, then $\IP(\mathfrak g)$ is
the flat limit of the corresponding $\IP^r$ in the fiber of the
classical $d$th Abel map on the generic (smooth) curve. Finally, we observe 
in Proposition \ref{prop:no-grds} that if a fiber of the $d$th Abel map 
for $X$ has any component of dimension less than $r$, then there is no limit 
linear series of dimension $r$ for the corresponding line bundle.

Finally, we mention that 
there 
is a somewhat parallel construction of
Eisenbud and Harris in Section 5 of \cite{e-h1}, where they describe
how the target projective space of the morphism associated to a linear
series degenerates in the special case that the limit is a refined limit
series. The space they describe is quite similar to the corresponding
special case of ours (see Remark \ref{rem:refined-case} below). However,
the constructions ought to be viewed
as dual to one another: if $\mathfrak g$ is a linear series on the generic
fiber, we construct degenerations of $\IP(\mathfrak g)$,
while the natural target of the associated morphism is not 
$\IP(\mathfrak g)$, but rather its dual space $\IP(\mathfrak g)^*$.

\section{Limit linear series}

Throughout this article, $X$ will denote the union of two smooth curves 
$Y$ and $Z$, meeting transversally at a point $P$. 

Let $\L$ be an invertible sheaf on $X$. It is determined by 
its restrictions $\L|_Y$ and $\L|_Z$. Also, there are natural short 
exact sequences,
\begin{equation}\label{LLY}
0\to\L|_Z(-P)\to\L\to\L|_Y\to 0,
\end{equation}\begin{equation}\label{LLZ}
0\to\L|_Y(-P)\to\L\to\L|_Z\to 0.
\end{equation}

For each integer $i$, let 
$\L^i$ be the invertible sheaf on $X$ with restrictions $\L|_Y(-iP)$ and 
$\L|_Z(iP)$. There are natural maps $\varphi^i\:\L^i\to\L^{i+1}$ and 
$\varphi_i\:\L^{i+1}\to\L^i$, defined as the compositions:
\begin{align*}
\varphi^i\:&\L^i\lra\L^i|_Z=\L|_Z(iP)=\L|_Z((i+1)P)(-P)=\L^{i+1}|_Z(-P)
\lra \L^{i+1},\\
\varphi_i\:&\L^{i+1}\lra\L^{i+1}|_Y=\L|_Y(-(i+1)P)=\L|_Y(-iP)(-P)=
\L^i|_Y(-P)\lra \L^i,
\end{align*}
where the first map in each composition is the restriction map, and the 

last maps are the inclusions in \eqref{LLY} and \eqref{LLZ} 
for $\L^{i+1}$ and $\L^i$ instead of $\L$. 
Notice that the compositions $\varphi^i\varphi_i$ and 
$\varphi_i\varphi^i$ are zero.

\begin{defn} 
Fix integers $d$ and $r$. A \emph{limit (linear) series} on $X$ of 
degree $d$ and dimension $r$ is a collection consisting of an invertible 
sheaf $\L$ on $X$ of degree $d$ on $Y$ and degree 0 on $Z$, and vector 
subspaces $V_i\subseteq\Gamma(X,\L^i)$ of dimension $r+1$, for each 
$i=0,\dots,d$, such that 
$\varphi^i(V_i)\subseteq V_{i+1}$ and $\varphi_i(V_{i+1})\subseteq V_i$ 
for each $i$.
\end{defn}

Given a limit series $(\L,V_0,\dots,V_d)$, we denote by $V_i^{Y,0}$ the 
subspace of $V_i$ of sections that vanish on $Y$, and by $V_i^{Z,0}$ the 
subspace of $V_i$ of sections that vanish on $Z$. 
Also, let $V_i|_Y$ denote the 
subspace of $\Gamma(Y,\L^i|_Y)$ generated by $V_i$ and $V_i|_Z$ that of 
$\Gamma(Z,\L^i|_Z)$ generated by the same $V_i$. 
Of course, $V_i^{Y,0}$ is the 
kernel of the surjection $V_i\to V_i|_Y$, and $V_i^{Z,0}$ is the 
kernel of the surjection $V_i\to V_i|_Z$. 
Also, the map $\varphi^i\: V_i\to V_{i+1}$ 
has kernel $V_i^{Z,0}$ and image contained in $V_{i+1}^{Y,0}$, whereas 
$\varphi_i\: V_{i+1}\to V_i$ has kernel $V_{i+1}^{Y,0}$ and 
image contained in $V_i^{Z,0}$. 

\begin{defn}\label{defn:exact} A limit linear series $(\L,V_0,\dots,V_d)$ 
is called \emph{exact} if, for each $i$, 
\begin{align*}
\Im(\vf^i\:V_i\to V_{i+1})=V_{i+1}^{Y,0}=\Ker(\vf_i\:V_{i+1}\to V_i),\\
\Im(\vf_i\:V_{i+1}\to V_i)=V_i^{Z,0}=\Ker(\vf^i\:V_i\to V_{i+1}).
\end{align*}
\end{defn}

It is a theorem of Liu \cite{li2}
that if $X$ is general (i.e., if both 
$(Y,P)$ and $(Z,P)$ are
general $1$-marked curves) the exact limit linear series are dense in
the space of all limit linear series.

One of the key properties of exact 
limit series
is the following, which may be
thought of as a simultaneous diagonalization lemma.

\begin{lemma}\label{lem:exact} If $(\L,V_0,\dots,V_d)$ is 
an exact limit series, 
then
there exist nonnegative integers $i_0 \leq i_1 \leq \dots \leq i_r \leq d$ 
and sections $s_0,\dots,s_r$ with $s_j \in V_{i_j}$ such that for each
$i=0,\dots,d$, the $s_j$ with $i_j=i$ form a basis of 
$V_i/(V_i^{Y,0} \oplus V_i^{Z,0})$, and the iterated images of all the 
$s_j$ form a basis for $V_i$.
\end{lemma}

For the argument, see the proof of Lemma A.12 (ii) of \cite{os8}.

\section{Abel maps}

To our knowledge, higher-degree Abel maps 
for curves of compact type appeared first in 
\cite{c-p1}, though they are the natural offspring of the 
construction of degree-1 Abel maps for stable curves in 
\cite{c-c-e1} or \cite{c-e1}.

We will need them in a very special situation, where they are 
easy to describe. 
Recall that $X$ is the union of two smooth curves $Y$ and $Z$, 
meeting transversally at a point $P$. Let $S^d(X)$ denote the symmetric 
product of $X$, thus parameterizing 0-cycles, or Weil divisors, on $X$ 
of degree $d$. 
The \emph{degree-$d$ Abel map} is a map 
$$
A_d\: S^d(X) \lra \Pic^d(X),
$$
where $\Pic^d(X)$ is the Picard scheme of $X$, parameterizing line 
bundles of a fixed multidegree $(d_1,d_2)$ with total degree 
$d_1+d_2=d$. The specific multidegree varies 
according to choices of components 
and polarizations; see \cite{c-p1}.

For our purposes, 
it is better to think of 
$\Pic^d(X)$ as parameterizing equivalence classes of line 
bundles of total degree $d$, where two line bundles 
$\L_1$ and $\L_2$ are said to be \emph{equivalent} 
if there exists an integer $j$ such that $\L_1|_Y\cong\L_2|_Y(-jP)$ and 
$\L_1|_Z\cong\L_2|_Z(jP)$. The map $A_d$ 
is then given as follows: Given a 0-cycle 
$D$ on $X$ of degree $d$, write it as $D=D_Y+D_Z$, 
where $D_Y$ and $D_Z$ are 
0-cycles, the first supported on $Y$ and the second on $Z$; then 
the image of 
$D$ under $A_d$ is the 
(class of the) 
line bundle on $X$ whose restrictions 
to $Y$ and $Z$ 
are $\O_Y(D_Y)$ and $\O_Z(D_Z)$. (Note the abuse of notation, where 
we view a 0-cycle of $X$ supported on $Y$ or $Z$ as a 0-cycle on $Y$ or 
$Z$, and vice versa.) This description of $A_d(D)$ 
does not depend on how $D$ is decomposed as $D_Y+D_Z$, 
as the line bundles 
resulting from different decompositions are all equivalent to each other. 

The fibers of $A_d$ are also easy to 
describe, at least set-theoretically. 
A point of $\Pic^d(X)$ 
has a unique representative $\L$ of degree $d$ on $Y$ and 0 on $Z$. 
Define the bundles $\L^i$ as in the last section. Set 
$\Gamma^i_Y:=\Gamma(Y,\L^{d-i}|_Y)$ and 
$\Gamma^i_Z:=\Gamma(Z,\L^i|_Z)$ for 
each $i=0,\dots,d$. There are natural closed embeddings 
$$
\IP(\Gamma^i_Y)\lra S^i(Y)\quad\text{and}\quad
\IP(\Gamma^i_Z)\lra S^i(Z),
$$
sending the class of a nonzero section $s$ to the 0-cycle $\div(s)$ 
associated to its zero scheme. Taking 
products
of these embeddings, 
and composing with the natural embeddings
$$
S^{d-i}(Y)\times S^i(Z)\lra S^d(X),
$$
which send a pair of 0-cycles to their sum, 
we obtain 
as images
subsets 
of $S^d(X)$, whose union is the fiber of 
$A_d$ over $\L$. Abusing notation, by not keeping record of the 
embeddings, we have:
\begin{equation}\label{eq:abel-fiber}
A_d^{-1}(\L)=\IP(\Gamma^0_Y)\times\IP(\Gamma^d_Z)\,\cup\,
\IP(\Gamma^1_Y)\times\IP(\Gamma^{d-1}_Z)\,\cup\cdots\cup\,
\IP(\Gamma^d_Y)\times\IP(\Gamma^0_Z).
\end{equation}

The natural inclusions 
$$
S^i(Y)\lra S^{i+1}(Y)\quad\text{and}\quad
S^i(Z)\lra S^{i+1}(Z),
$$
sending a 0-cycle $D$ to $D+P$ in both cases, 
take $\IP(\Gamma^i_Y)$ to $\IP(\Gamma^{i+1}_Y)$ 
and $\IP(\Gamma^i_Z)$ to $\IP(\Gamma^{i+1}_Z)$, respectively. 
Abusing notation again, inside 
$S^d(Y)$ and $S^d(Z)$ we have chains of subschemes:
\begin{equation}\label{PY}
\IP(\Gamma^0_Y)\subseteq\IP(\Gamma^1_Y)\subseteq\cdots\subseteq
\IP(\Gamma^{d-1}_Y)\subseteq\IP(\Gamma^d_Y),
\end{equation}
\begin{equation}\label{PZ}
\IP(\Gamma^0_Z)\subseteq\IP(\Gamma^1_Z)\subseteq\cdots\subseteq
\IP(\Gamma^{d-1}_Z)\subseteq\IP(\Gamma^d_Z).
\end{equation}

So we may consider the union on the right-hand side of 
\eqref{eq:abel-fiber} inside the product 
$\IP(\Gamma^d_Y)\times\IP(\Gamma^d_Z)$ instead of $S^d(X)$. 
It is equal to $A_d^{-1}(\L)$ nonetheless. Indeed, 
the product $S^d(Y)\times S^d(Z)$ can be viewed 
naturally inside 
$S^{2d}(X)$, by sending a pair of 0-cycles to their sum. 
Also, $S^d(X)$ can be viewed inside $S^{2d}(X)$, by 
sending a 0-cycle $D$ to $D+dP$. 
Under these inclusions, 
each $\IP(\Gamma^i_Y)\times\IP(\Gamma^{d-i}_Z)$, 
whether viewed as a subset of $S^d(X)$ or of 
$\IP(\Gamma^d_Y)\times\IP(\Gamma^d_Z)$, gives the same subset of 
$S^{2d}(X)$.

Given 
$R\in X-P$, 
let
$$
R+S^{d-1}(X):=\{D\in S^d(X)\,|\, D\geq R\}.
$$
Fixing points 
$R_1\in Y-P$ and $R_2\in Z-P$, 
consider the 
divisors $H_{s,t}$ on $S^d(X)$ given by
$$
H_{s,t}:=s(R_1+S^{d-1}(X))+t(R_2+S^{d-1}(X))
$$ 
for each $s,t\in\mathbb Z$. 
If $s,t>0$, the restriction of $H_{s,t}$ to any 
fiber of $A_d$ is ample.

Given a subscheme 
$W$ of 
a fiber $A_d^{-1}(\L)$, 
we may compute its 
\emph{bivariate Hilbert polynomial}
$$
P_W(s,t):=\dim\Gamma(W,\O_W(s,t))=\dim\Gamma(W,\O_{S^d(X)}(H_{s,t})|_W)
$$
for $s,t>>0$. 
Notice that, viewing $W$ inside $\IP(\Gamma^d_Y)\times\IP(\Gamma^d_Z)$ 
under the 
natural embedding 
$A_d^{-1}(\L)\to\IP(\Gamma^d_Y)\times\IP(\Gamma^d_Z)$,
the Hilbert 
polynomial $P_W$ is the bivariate Hilbert polynomial of a subscheme 
in the product of two projective spaces, and is thus independent of 
the choices of $R_1$ and $R_2$. 

\section{Limit linear series and fibers of the Abel map}

\begin{defn} Let $\mathfrak g:=(\L,V_0,\dots,V_d)$ 
be an exact limit linear series on $X$ 
of degree $d$ and dimension $r$. Let $\IP(\mathfrak g)$ denote the 
closure of the subset of $A_d^{-1}(\L)$ consisting of points of the form
$$
\div(s|_Y)+\div(s|_Z)\in S^d(X),
$$
for $s\in V_i-(V_i^{Y,0}\cup V_i^{Z,0})$ for some $i$. We give 
$\IP(\mathfrak g)$ the reduced induced subscheme structure.
\end{defn}

The definition makes sense, because, for each $i$, the line bundle 
associated to $\div(s|_Y)+\div(s|_Z)$ under $A_d$ is clearly $\L^i$, 
which is equivalent to $\L$. 

\begin{rmk}\label{rem:empty-case}
Note that the definition of $\IP(\mathfrak g)$ will only yield divisors
from $V_i$ if $V_i \neq V_i^{Y,0} \cup V_i^{Z,0}$. However, this condition 
is weak and easy to control: 
it is violated if and only if 
$V_i=V_i^{Y,0}$ or
$V_i=V_i^{Z,0}$, or equivalently, 
if and only if 
$V_i|_Z=0$ or $V_i|_Y=0$. Moreover, 
since the $V_i$ give filtrations of $V_0|_Y$ and $V_d|_Z$, 
violation occurs only 
in ``extremal'' degrees. 
More precisely, there exist certain
$i_{\low} \leq i_{\high}$ between $0$ and $d$ such that 
$V_i \neq V_i^{Y,0} \cup V_i^{Z,0}$ if and only if 
$i_{\low} \leq i \leq i_{\high}$.

Furthermore, exactness yields 
$$
V_{i_{\low}} \neq V_{i_{\low}}^{Y,0} \oplus V_{i_{\low}}^{Z,0}
\quad\text{and}\quad
V_{i_{\high}} \neq V_{i_{\high}}^{Y,0} \oplus V_{i_{\high}}^{Z,0}.
$$
Indeed, if
$i_{\low}=0$, then $V_0^{Y,0}=0$.
Otherwise, since $V_{i_{\low}-1}^{Z,0}=V_{i_{\low}-1}$, the image of 
$V_{i_{\low}-1}$ in $V_{i_{\low}}$ is $0$,
whence 
by exactness $V_{i_{\low}}^{Y,0}=0$. In either case, since
$V_{i_{\low}}^{Z,0} \neq V_{i_{\low}}$ by 
definition of $i_{\low}$,
we get the desired
assertion for $V_{i_{\low}}$. The same argument works for $V_{i_{\high}}$,
switching the role of $Y$ and $Z$.
\end{rmk}

\begin{thm}\label{main} If $\mathfrak g=(\L,V_0,\dots,V_d)$ 
is an exact limit linear series on $X$ 
of degree $d$ and dimension $r$, then $\IP(\mathfrak g)$ is reduced, connected,
and Cohen--Macaulay of 
pure dimension $r$, has bivariate 
Hilbert polynomial $P(s,t)=\binom{s+t+r}{r}$ 
and is a flat degeneration of $\IP^r$.
\end{thm}

The proof of the theorem is lengthy, 
so
we break it into a number of
steps. The first three lemmas are independent of limit linear series,
and will ultimately be used to describe the geometry of 
$\IP(\mathfrak g)$, starting with its irreducible components.

\begin{lemma}\label{lem:Qi} Given $p,q,m$ nonnegative integers, 
let $Q_{p,q,m} \subseteq \IP^{m+q} \times \IP^{m+p}$ be the subscheme 
given by 
the equations 
$$
x_i y_j - x_j y_i = 0 \quad\text{for }i,j=p,\dots,m+p,
$$
where $\IP^{m+q}$ is given coordinates $x_{p},\dots,x_{m+p+q}$, and
$\IP^{m+p}$ is given coordinates $y_0,\dots,y_{m+p}$.
Then $Q_{p,q,m}$ is integral 
of dimension $m+p+q$, 
determinantal 
and hence Cohen--Macaulay, with bivariate Hilbert polynomial 
$P_{Q_{p,q,m}}(s,t)$ given by
$$\sum_{\ell=0}^m \binom{s+q+\ell}{q+\ell}\binom{t+p+m-\ell}{p+m-\ell}
-\sum_{\ell=0}^{m-1}
\binom{s+q+\ell}{q+\ell}\binom{t+p+m-1-\ell}{p+m-1-\ell}.$$
\end{lemma}

(The variable indexing in the lemma may appear \textit{ad hoc}, but it
will be useful when we consider certain unions of these varieties inside
larger products of projective spaces.)

\begin{proof} Denote by $S$ the bigraded ring 
$k[x_{p},\dots,x_{m+p},y_{p},\dots,y_{m+p}]/I$, where $I$ is
the ideal generated by
$$
x_i y_j - x_j y_i = 0 \quad\text{for }i,j=p,\dots,m+p.
$$
Then $S$ is the homogeneous coordinate ring of the diagonal in
$\IP^{m} \times \IP^{m}$, and its bivariate Hilbert polynomial is
$\binom{s+t+m}{m}$. Moreover, since the cohomology of projective space
vanishes in positive degree for $\O(s+t)$ with $s+t \geq 0$, we have that
the Hilbert polynomial agrees with the Hilbert function.  

The homogeneous coordinate ring of $Q_{p,q,m}$ is then
$$S_{Q}:=S[x_{m+p+1},\dots,x_{m+p+q},y_{0},\dots,y_{p-1}],$$
from which it follows that $Q_{p,q,m}$ is integral. Furthermore, 
a monomial of bidegree
$(s,t)$ in $S_Q$ may be written uniquely as a monomial of bidegree $(i,j)$
in $S$ times a monomial of degree $s-i$ in the 
$x_{n}$ for $n>m+p$,
and a monomial of degree $t-j$ in the 
$y_{n}$ for $n<p$, 
where
$i$ and $j$
are 
nonnegative integers less than or equal to $s$ and $t$
respectively. But we know that there are $\binom{i+j+m}{m}$ monomials of 
bidegree $(i,j)$ in $S$, and this is equal to the number of monomials of 
degree $i+j$ in $m+1$ variables,
say $z_0,\dots,z_m$. We thus conclude that the number of monomials of
bidegree $(s,t)$ in $S_Q$ is equal to the number of monomials of total
degree $s+t$ in 
$x_{m+p+1},\dots,x_{m+p+q},y_0,\dots,y_{p-1},z_0,\dots,z_m$, 
with the degrees in the 
$x_{n}$
and $y_{n}$ bounded by $s$ and $t$, respectively. Denote the latter set
of monomials by $\mathcal S_Q(s,t)$. In the definition of $\mathcal S_Q(s,t)$ 
and hereafter the $x_{n}$ are assumed to have $n>m+p$ and the $y_n$ are 
assumed to have $n<p$.

On the other hand, each summand in the first sum of the desired formula
for the Hilbert polynomial is equal to the product of the number of 
monomials of degree $s$ in the $x_n$ 
and 
$z_0,\dots,z_{\ell}$ 
with the number of monomials of degree $t$ in
the 
$y_n$ and $z_{\ell},\dots,z_m$. 
Taking the product of such a pair of
monomials trivially gives an element of $\mathcal S_Q(s,t)$. For a given 
$\ell$, the resulting map from pairs of monomials to $\mathcal S_Q(s,t)$ is 
clearly 
injective. And if we take the union of all these maps over all $\ell$,
it is clear that 
we get a surjective map, 
as for each monomial in $\mathcal S_Q(s,t)$ we may choose 
$\ell$ minimal so that 
its 
total degree 
in the $x_n$ and $z_0,\dots,z_{\ell}$ is at least $s$. 

However, a given monomial in $\mathcal S_Q(s,t)$ may arise in this way from 
more than one pair of monomials. The number of times 
it
arises is 
equal to the number of $\ell$ such that the total degree in the $x_n$
and $z_0,\dots,z_{\ell}$ is at least $s$ and the total degree in the
$y_n$ and $z_{\ell},\dots,z_m$ is at least $t$. If this holds for 
$\ell',\ell'+1,\dots,\ell''$, with $\ell'<\ell''$, it immediately follows 
that the 
total degree of the monomial in the $x_n$ and $z_0,\dots,z_{\ell'}$ is 
$s$, its total degree in $z_{\ell'+1},\dots,z_{\ell''-1}$ is 0, and 
its total degree in the $y_n$ and $z_{\ell''},\dots,z_m$ is $t$. 
Each summand 
in the second sum of the desired formula may be interpreted as counting
pairs of monomials of degree $s$ in the $x_n$ and $z_0,\dots,z_{\ell}$
and of degree $t$ in the $y_n$ and $z_{\ell+1},\dots,z_m$. These
likewise map to $\mathcal S_Q(s,t)$, and we see that the number mapping to a
given monomial is precisely $\ell''-\ell'$, with notation as above. 
We thus conclude the desired formula for the
Hilbert polynomial. 
And it follows from the formula that the dimension of $Q_{p,q,m}$ is 
the one stated.

Now, the description of 
$Q_{p,q,m}$ is visibly determinantal,
coming from the condition that the matrix
$$
\left[\begin{matrix}
x_{p}&x_{p+1}&\dots&x_{m+p-1}&x_{m+p}\\
y_{p}&x_{p+1}&\dots&y_{m+p-1}&y_{m+p}
\end{matrix}\right]
$$
have rank at most $1$. Moreover,
$$
\text{codim}(Q_{p,q,m},\IP^{m+q}\times\IP^{m+p})
=(m+q)+(m+p)-(m+p+q)= m,
$$
and the expected codimension of $Q_{p,q,m}$ is $(2-1)((m+1)-1)=m$, 
so we conclude that $Q_{p,q,m}$ is determinantal 
of the expected codimension, 
and consequently
Cohen--Macaulay.
\end{proof}

\begin{lemma}\label{lem:union} 
Given nonnegative integers $r$ and $m_0,\dots,m_n$ 
with $\sum_i (m_i+1) \leq r+1$, define sequences $p_0,\dots,p_n$ and 
$q_0,\dots,q_n$ as follows: 
\begin{enumerate}
\item $p_0:=0$ and $p_i:=p_{i-1}+m_{i-1}+1$ for each $i=1,\dots,n$;
\item $q_i:=r-m_i-p_i$ for each $i=0,\dots,n$.
\end{enumerate}
Give $\IP^r \times \IP^r$ 
bihomogeneous 
coordinates $x_j$ and $y_j$ for $j=0,\dots,r$,
and for each $i=0,\dots,n$, 
view
the $Q_{p_i,q_i,m_i}$ of Lemma
\ref{lem:Qi} as a subvariety $Q_i$ of $\IP^r \times \IP^r$ in the 
natural
way, by considering 
$$\IP^{m_i+q_i}\times \IP^{m_i+p_i}
=V(x_0,\dots,x_{p_i-1},y_{m_i+p_i+1},\dots,y_r) 
\subseteq \IP^r \times \IP^r.$$
Set 
$$Q = 
\bigcup_{i=0}^n 
Q_i \subseteq \IP^r \times \IP^r.$$
Then $Q$ is connected, reduced and Cohen--Macaulay
of dimension $r$, with bivariate Hilbert polynomial $P_Q(s,t)$ given by
$$\sum_{\ell=0}^{p_n+m_n}\binom{s+r-\ell}{r-\ell}
\binom{t+\ell}{\ell}
-\sum_{\ell=0}^{p_n+m_n-1}
\binom{s+r-1-\ell}{r-1-\ell}\binom{t+\ell}{\ell}.$$
If further $p_n+m_n=r$, 
then
$$P_Q(s,t)=\binom{s+t+r}{r}.$$
\end{lemma}

(Note that $p_n+m_n+1=\sum_i (m_i+1)$, so by hypothesis the case $p_n+m_n=r$
is maximal.)

\begin{proof} We begin by observing that the final assertion, for the case 
$p_n+m_n=r$, follows immediately from the rest of the lemma. Indeed, 
by
the general formula, our Hilbert polynomial in this case
agrees with that 
in the case $n=0,m_0=r$, in which $Q=Q_0$ is simply
the diagonal in $\IP^r \times \IP^r$, and therefore has the desired Hilbert 
polynomial.

For the main statement of the lemma,
the case $n=0$ follows immediately from Lemma \ref{lem:Qi}. 
(The stated form of the Hilbert polynomial follows from that in 
Lemma~\ref{lem:Qi} by replacing $\ell$ by $m-\ell$ in the first sum and 
$\ell$ by $m-1-\ell$ in the second.)
Now, suppose 
$n>0$. Denote by $Q'$ the union of the $Q_{i}$ for $i=0,\dots,n-1$,
so that $Q=Q' \cup Q_n$. 
Set $Q'':=Q' \cap Q_n$. We claim that (scheme-theoretically)
$$Q''=V(x_0,\dots,x_{p_n-1},y_{p_n},\dots,y_r)\cong
\IP^{r-p_{n}} \times \IP^{p_{n}-1}.$$ 
The lemma follows from the claim and induction on $n$. Indeed, 
$Q_n$ is reduced of pure dimension $r$ by Lemma \ref{lem:Qi}. So is 
$Q'$ by induction, 
and thus so is their union $Q$. 
Also, $Q'$ and $Q_n$ 
are connected, and thus so is $Q$ because $Q''$ is nonempty by the claim. 
Furthermore, $Q'$ and $Q_n$ are Cohen--Macaulay, and thus so is $Q$ by 
the claim and 
\cite{ei1}, Ex.~18.13, p.~467.
Finally,
we 
obtain the desired form for the Hilbert polynomial simply by
applying 
induction, the claim,
and the identity
$$P_Q(s,t)=P_{Q'}(s,t)+P_{Q_n}(s,t)-P_{Q''}(s,t).$$

We now prove the claim.
Denote by $J_i$ the ideal defining $Q_i$, so that
$$
J_i=(x_0,\dots,x_{p_i-1},y_{m_i+p_i+1},\dots,y_r)
+(x_j y_\ell - x_j y_\ell \,|\,p_i\leq j,\ell\leq m_i+p_i).
$$
We need 
only
check that
$$
J_n+(J_{n-1}\cap\cdots\cap J_0)=J_n+J_{n-1}=
(x_0,\dots,x_{p_{n}-1},y_{m_{n-1}+p_{n-1}+1},\dots,y_r),
$$
from which the claim follows by noting that 
$p_n=m_{n-1}+p_{n-1}+1$.
The second equality above is trivial. As for the first, we
clearly have
$J_n+(J_{n-1}\cap\cdots\cap J_0)\subseteq J_n+J_{n-1}$. 
To show the reverse inclusion we observe that
$x_0,\dots,x_{p_n -1}\in J_n$, while 
$y_{m_{n-1}+p_{n-1}+1},\dots,y_r\in J_{n-1}\cap\cdots\cap J_0$.
Thus we have proved the claim and the lemma.
\end{proof}

\begin{lemma}\label{lem:degeneration} In the situation of Lemma
\ref{lem:union}, if $p_n+m_n=r$ then $Q$ is a flat degeneration of the
diagonal in $\IP^r \times \IP^r$.
\end{lemma}

\begin{proof}
Consider the subscheme $\widehat{W}$ of 
$\IP^r \times \IP^r \times (\mathbb A^1_z\smallsetminus(0))$ given 
by the equations 
$z^{\epsilon_j} x_iy_j-z^{\epsilon_i}x_jy_i=0$ for $0\leq i\leq j\leq r$,
where 
$\epsilon_j=n-i$ if $p_i\leq j \leq p_i+m_i$.
At $z=1$, we recover the equations of the diagonal in $\IP^r \times \IP^r$,
and we see moreover that for any $z \neq 0$, the scheme we obtain is simply 
a change of coordinates of the diagonal, given by 
$y_j \mapsto z^{\epsilon_j}y_j$. Let $W$ be the closure of $\widehat{W}$ in
$\IP^r \times \IP^r \times \mathbb A^1_z$, and $W_0$ the fiber over $z=0$.
Then $W_0$ is a flat degeneration of the diagonal. Thus, by Lemma 
\ref{lem:union} we have that 
$W_0$ has the same Hilbert polynomial as $Q$. 
Since, by the same lemma, $Q$ is reduced,
if we show that $Q \subseteq W_0$ set-theoretically,
we conclude scheme-theoretic equality, and the lemma
is proved.

Now, by definition a point $((a_0,\dots,a_r),(b_0,\dots,b_r))$ of $Q$ lies 
in $Q_i$ for some $i$. We then have 
$$a_0=\dots=a_{p_i-1}=b_{m_i+p_i+1}=\dots=b_r=0;$$
also, 
there exist $\lambda, \mu$ not both zero with
$\lambda a_i = \mu b_i$ for all $i=p_i,\dots,m_i+p_i$. Since $Q_i$ is 
irreducible by Lemma \ref{lem:Qi}, it is enough to show that a nonempty
open subset of $Q_i$ lies in $W_0$, so we may further assume that 
neither $(a_{p_i},\dots,a_{m_i+p_i})$ nor $(b_{p_i},\dots,b_{m_i+p_i})$ is identically zero, 
so that $\lambda$ and $\mu$ are both nonzero. Then we
can renormalize the $a_i$ and $b_i$ so that $\lambda=\mu=1$.
In order to show that our point is in $W_0$, we define $e_j=-\epsilon_j$
for $j \geq p_i$, and $e_j=-\epsilon_{p_i}$ for $j<p_i$. Then set
$$c_j
=\begin{cases} a_j z^{e_j}:& j \geq p_i \\ b_j z^{e_j}:&\text{otherwise.}
\end{cases}$$
Then we evidently have the section
$$((c_0 z^{\epsilon_0},\dots,c_r z^{\epsilon_r}),(c_0,\dots,c_r))$$
in $W$ for all $z \neq 0$, and 
hence its limit at $z=0$ is 
in $W_0$ by definition. To obtain the limit at $z=0$ we have to
divide through on the left and right by the minimum powers of $z$ occurring
on each side. On the right, since the $\epsilon_j$ are 
nonincreasing 
in $j$,
this minimum is $-\epsilon_{p_i}$, achieved precisely for 
the $j\in\{0,\dots,p_i+m_i\}$ such that $b_j\neq 0$.
We thus see that the limit on the right is 
$$(b_0,\dots,b_{p_i-1},a_{p_i},\dots,a_{p_i+m_i},0\dots0)=(b_0,\dots,b_r).$$
Similarly, the minimum on the left is $0$, achieved precisely for 
the $j\in\{p_i,\dots,r\}$ such that $a_j\neq 0$, 
so the limit is
$$(0,\dots,0,a_{p_i},\dots,a_r)=(a_0,\dots,a_r).$$
We have thus shown that our chosen point is in $W_0$, and we conclude
that $Q_i \subseteq W_0$ for each $i$, and thus that $Q \subseteq W_0$,
as desired.
\end{proof}

We now relate the previous lemmas to $\IP(\mathfrak g)$. 
Given $\mathfrak g$, we observe that for any $i$ 
we have a closed embedding 
\begin{equation}\label{fiber-prod} \IP(V_i|_Y) \times \IP(V_i|_Z) \lra
A_d^{-1}(\L);\end{equation}
see \eqref{eq:abel-fiber}.

\begin{defn} Let 
$\mathfrak g=(\L,V_0,\dots,V_d)$ 
be an exact limit
linear 
series. For each $i$ such that $V_i \neq V_i^{Y,0} \cup V_i^{Z,0}$, denote 
by $\IP(\mathfrak g_i)$ the 
closure of the subset of $A_d^{-1}(\L)$ consisting of points of the form
$$
\div(s|_Y)+\div(s|_Z)\in S^d(X),
$$
for $s\in V_i-(V_i^{Y,0}\cup V_i^{Z,0})$.
\end{defn}

Observe that the construction of $\IP(\mathfrak g_i)$ visibly factors through
\eqref{fiber-prod}, so we have
$$\IP(\mathfrak g_i) \subseteq \IP(V_i|_Y) \times \IP(V_i|_Z).$$

\begin{lemma}\label{lem:get-Qi} Let 
$\mathfrak g=(\L,V_0,\dots,V_d)$ 
be an exact limit
linear 
series. Let $i$ be 
such that
$V_i \neq V_i^{Y,0} \cup V_i^{Z,0}$.
Set $p:=\dim V_i^{Y,0}$, $q:=\dim V_i^{Z,0}$ and $m:=r-p-q$, and let
$s_0,\dots,s_r \in V_i$ be any basis such that $s_0,\dots,s_{p-1}$ form
a basis for $V_i^{Y,0}$, and $s_{r+1-q},\dots,s_r$ form a basis for
$V_i^{Z,0}$. Let the $s_i$ induce coordinates $x_p,\dots,x_r$
on $\IP(V_i|_Y)$ and $y_0,\dots,y_{r-q}$ on $\IP(V_i|_Z)$.

Then if $m=-1$, we have 
$$\IP(\mathfrak g_i) = \IP(V_i|_Y) \times \IP(V_i|_Z).$$ 
If $m \geq 0$, we have 
$$\IP(\mathfrak g_i)=Q_{p,q,m} \subseteq \IP(V_i|_Y) \times \IP(V_i|_Z),$$ 
with notation as in Lemma \ref{lem:Qi}.
\end{lemma}

(Note that since $V_i^{Y,0} \cap V_i^{Z,0}=(0)$, we necessarily have
$p+q \leq r+1$, with equality occurring when $\IP(V_i|_Y) \times \IP(V_i|_Z)$
has dimension $r-1$.)

\begin{proof} 
Since the kernel of $V_i \to V_i|_Y$ is by definition $V_i^{Y,0}$, we may 
view the points of 
$\IP(V_i|_Y)$ as the subspaces of $V_i$ containing $V_i^{Y,0}$ with 
codimension $1$. An analogous interpretation holds 
for $\IP(V_i|_Z)$. With this point of view, consider the subset 
$$
Q_i\subseteq \IP(V_i|_Y)\times\IP(V_i|_Z)
$$
parameterizing pairs $(W_1,W_2)$ of subspaces of $V_i$, with 
$V_i^{Y,0}\subseteq W_1$ and $V_i^{Z,0}\subseteq W_2$ of codimension $1$, 
such that $W_1\cap W_2$ is nontrivial. Then $Q_i$ is a closed subset.
If $p+q=r+1$, it is clear that 
$Q_i=\IP(V_i|_Y) \times \IP(V_i|_Z)$;  
otherwise we claim that $Q_i$
agrees with the $Q_{p,q,m}$ of Lemma \ref{lem:Qi}. 

Indeed, in 
the coordinates $x_j$, we may represent a point of 
$\IP(V_i|_Y)$ by $(a_p,\dots,a_r)$. The corresponding subspace of $V_i$
is then the span of $s_0,\dots,s_{p-1},\sum_{i \geq p} a_i s_i$.
Similarly, 
the point
$(b_0,\dots,b_{r-q}) \in \IP(V_i|_Z)$ corresponds to the 
subspace of $V_i$ spanned by
$\sum_{i \leq r-q} b_i s_i, s_{r-q+1},\dots, s_r$.
The intersection of these subspaces is nontrivial if and only if 
the span of all the above sections generate a subspace 
of dimension at most $p+q+1$, thus if and only if the 
matrix
$$
\left[\begin{matrix}
a_{p}&a_{p+1}&\dots&a_{r-q-1}&a_{r-q}\\
b_{p}&b_{p+1}&\dots&b_{r-q-1}&b_{r-q}
\end{matrix}\right]
$$
has rank at most $1$. Thus, in the coordinates $x_j, y_m$,
we see that $Q_i$ is the closed subset 
given by the equations
$$
x_jy_m-x_my_j=0\quad\text{for }j,m=p,\dots,r-q.
$$
We thus conclude that $Q_i=Q_{p,q,m}$, as claimed. In particular, by
Lemma \ref{lem:Qi}, we have that $Q_i$ is integral.

It thus suffices to show (set-theoretically) that $\IP(\mathfrak g_i)=Q_i$. 
One containment is clear: let
$v\in V_i-(V_i^{Y,0}\cup V_i^{Z,0})$. Then $W_1:=kv+V_i^{Y,0}$ and 
$W_2:=kv+V_i^{Z,0}$ are subspaces of $V_i$ with nontrivial intersection, 
and thus the pair $(W_1,W_2)$ defines a point of $Q_i$. So 
$\IP(\mathfrak g_i) \subseteq Q_i$. Since $\IP(\mathfrak g_i)$ is defined 
via closure, it suffices to
prove that it contains a nonempty open subset of $Q_i$.
First, if $p+q=r+1$, we have $V_i=V_i^{Z,0} \oplus V_i^{Y,0}$, so
any $(W_1,W_2)\in Q_i=\IP(V_i|_Y) \times \IP(V_i|_Z)$ may be represented 
by $(w_1,w_2) \in 
(V_i^{Z,0}\smallsetminus (0)) \times (V_i^{Y,0} \smallsetminus (0)),$
and then $(W_1,W_2)$ is the image of $w_1+w_2$.
Thus, we see in this case that 
$Q_i$ coincides with the set whose closure defines $\IP(\mathfrak g_i)$, 
so closure is unnecessary in this case, and in particular 
$Q_i=\IP(\mathfrak g_i)$.
On the other hand, if $p+q < r+1$, then there is a nonempty open subset
of
$Q_i$ consisting of $(W_1,W_2)$ such that
$W_1\cap V_i^{Z,0}=0$ and $W_2\cap V_i^{Y,0}=0$. On this subset, any nonzero 
$v\in W_1\cap W_2$ is neither in $V_i^{Y,0}$ nor in $V_i^{Z,0}$, and 
thus $W_1=kv+V_i^{Y,0}$ and $W_2=kv+V_i^{Z,0}$. So 
$(W_1,W_2) \in \IP(\mathfrak g_i)$.
We thus conclude $\IP(\mathfrak g_i)=Q_i$, as desired.
\end{proof}

We are now ready to complete the proof of our main theorem.

\begin{proof}[Proof of Theorem \ref{main}] Let $i_0,\dots,i_r$
and $s_0,\dots,s_r$ be as in 
Lemma \ref{lem:exact}. So the $s_j$ 
yield bases for each $V_i$, and in particular for $V_0$ and $V_d$.
Since 
the restriction maps $V_0 \to V_0|_Y$ and $V_d \to V_d|_Z$ are isomorphisms, 
the $s_j$ induce coordinates on $\IP(V_0|_Y)$ and
$\IP(V_d|_Z)$, to be denoted $x_0,\dots,x_r$ and
$y_0,\dots,y_r$ respectively. 

Let $m_0$ be the number of $j>0$ with $i_j=i_0$, let $m_1$ be the number
of $j>m_0+1$ with $i_j=i_{m_0+1}$, and so forth. Suppose that $n+1$
is the number of distinct values of the $i_j$, so that we obtain a
sequence of nonnegative numbers $m_0,\dots,m_n$.
In light of Lemmas \ref{lem:union} and \ref{lem:degeneration}, the 
theorem will follow if we show that under the coordinates
$x_i$, $y_i$, we have that
$\IP(\mathfrak g)$ coincides with the $Q$ of Lemma \ref{lem:union}. 
Following through the definitions and applying Lemma \ref{lem:get-Qi},
we see that $Q$ is precisely the union of the $\IP(\mathfrak g_i)$ 
for those $V_i$ with $i=i_j$ for some $j$. It thus suffices 
to show that if $i \neq i_j$ for any $j$, we get nothing new from
$V_i$, that is, $\IP(\mathfrak g_i)$ is already contained in some 
$\IP(\mathfrak g_{i_j})$.

According to Remark \ref{rem:empty-case}, if we either increase or
decrease $i$ we will eventually get to some $i'$ with
$V_{i'} \neq V_{i'}^{Y,0} \oplus V_{i'}^{Z,0}$, or equivalently
$i'=i_j$ for some $j$. Suppose we chose to decrease 
$i$.
Thus
$i_j < i$, and there is no $i_{j'}$ with $i_j < i_{j'} \leq i$. It
follows 
from exactness
that for $i'>i_j$, with $i' \leq i$, we have $V_{i'}^{Y,0}$
identified with $V_i^{Y,0}$ under the map $V_{i'} \to V_i$, and 
$V_{i}^{Z,0}$ identified with $V_{i'}^{Z,0}$ under the map
$V_i \to V_{i'}$. Thus, each $\IP(\mathfrak g_{i'})$ is equal to
$\IP(\mathfrak g_i)$. Furthermore, we still have 
$V_i^{Z,0}$ identified with $V_{i_j}^{Z,0}$,
and the codimension of the image of 
$V_{i_j}^{Y,0}$ in $V_i^{Y,0}$ 
is equal to 
$\dim V_{i_j}/(V_{i_j}^{Y,0} \oplus V_{i_j}^{Z,0})=m_{\ell}$, 
if $i_j$ is the $(\ell+1)$st value taken on by the $i_{j'}$.
It then follows from the explicit 
descriptions given in Lemma \ref{lem:get-Qi} 
that 
$\IP(\mathfrak g_i) \subseteq \IP(\mathfrak g_{i_j})$.
We conclude that 
$\IP(\mathfrak g_i)\subseteq Q$, 
and thus the theorem.
\end{proof}

\begin{rmk}\label{rem:refined-case} It follows from the above proof that 
$$
\IP(\mathfrak g)=\bigcup_{\ell=0}^m\IP(\mathfrak g_{i_\ell}),
$$
where $i_0,\dots,i_m$ is the sequence of integers $i=0,\dots,d$ 
such that $V_i \neq V_i^{Y,0} \oplus V_i^{Z,0}$.

In the special case that $V_i^{Y,0} \oplus V_i^{Z,0}$ has 
codimension 1 in $V_i$, it follows from Lemma \ref{lem:get-Qi} 
that $\IP(\mathfrak g_i)$ is isomorphic to an $r$-dimensional product of 
two projective spaces. 

Such a situation arises frequently: 
Indeed, for a refined limit linear series 
(see Definition 6.5 of \cite{os8}) the codimension of 
$V_i^{Y,0} \oplus V_i^{Z,0}$ in $V_i$ is either 0 or 1 for every $i$. 
Thus, in this case, $\IP(\mathfrak g)$ consists of a union
of $r+1$ irreducible components, each isomorphic to an $r$-dimensional
product of two projective spaces. 
\end{rmk}

\section{Limits of linear series}\label{sec:limits}

Let $B$ be the spectrum of a discrete valuation ring with algebraically 
closed residue field, and let $\eta$ denote its generic point. 

\begin{defn} Let $\pi\:\X\to B$ be a flat, projective map, where 
$\X$ is regular, the generic fiber of $\pi$ is 
smooth, and the special fiber is isomorphic to $X$, the union of two 
smooth curves $Y$ and $Z$ meeting transversally at a point $P$. We 
call $\pi$ or $\X/B$ a \emph{regular smoothing} of $X$.
\end{defn}

Let $\L_\eta$ be an invertible sheaf of degree $d$ on the generic 
fiber $\X_\eta$. Since $\X$ is regular, $\L_\eta$ extends to an 
invertible sheaf $\L$ on $\X$. Since $\X$ is regular, $Y$ and $Z$ are 
Cartier divisors on $\X$, and thus $\L(-iZ):=\L\ox\O_\X(-iZ)$ are 
also extensions of 
$\L_\eta$ for all $i$.
Thus there is an extension $\L$ such that 
$\L|_X$ has degree $d$ on $Y$ and 0 on $Z$. Fix this extension $\L$, 
and set $\L^i:=\L(-iZ)|_X$ for $i=0,\dots,d$. 

Let $V_\eta$ be a vector subspace of $\Gamma(\X_\eta,\L_\eta)$ of 
dimension $r+1$. Viewing $V_\eta$ as a subspace of 
$\Gamma(\X_\eta,\L(-iZ)|_{\X_\eta})$ for each $i=0,\dots,d$, set 
$$
\wt V_i:=\Gamma(\X,\L(-iZ))\cap V_\eta,
$$
and denote by $V_i\subseteq\Gamma(X,\L^i)$ 
the image of the restriction of $\wt V_i$ to the 
special fiber. 

Notice that 
$$
\wt V_i\cap\Gamma(\X,\L(-(i+1)Z))=\wt V_{i+1}.
$$
Thus, not only is the image of $V_{i+1}$ under the map on global sections 
induced by the natural map $\L^{i+1}\to\L^i$ contained in $V_i$, but it is
also equal to the subspace of $V_i$ of sections that vanish on $Z$. 
An analogous statement can be made with regard to the natural map 
$\L^i\to\L^{i+1}$ in the reverse direction. So
$$
\mathfrak g:=(\L^0,V_0,V_1,\dots,V_d)
$$
is an exact limit linear series. We say that $\mathfrak g$ is the 
\emph{limit} of the linear series $(\L_\eta,V_\eta)$.

\begin{thm}\label{thm:flat-limit} Let $\X/B$ be a regular smoothing of $X$ and 
$(\L_\eta,V_\eta)$ a linear series of dimension $r$ and degree $d$ 
on the generic fiber. Let $\mathfrak g$ be the limit linear series 
that is limit of $(\L_\eta,V_\eta)$. Then $\IP(V_\eta)$, viewed as a 
subscheme of the fiber of the 
relative symmetric product $S^d(\X/B)$ over $\eta$, has 
closure 
intersecting $S^d(X)$ in
$\IP(\mathfrak g)$.
\end{thm}

\begin{proof} Up to making an \'etale base change, we may assume that 
$\X/B$ has two sections $\Sigma_1$ and $\Sigma_2$, the first intersecting 
$Y$ away from $P$, the second intersecting $Z$ away from $P$. Let 
$$
\H_{s,t}:=s(\Sigma_1+S^{d-1}(\X/B))+t(\Sigma_2+S^{d-1}(\X/B))
$$
for each $s,t\in\mathbb Z$. The restriction of $\H_{s,t}$ to 
$\IP(V_\eta)$ is equal to $\O_{\IP(V_\eta)}(s+t)$, thus
$$
P_{\IP(V_\eta)}(s,t):=
\dim\Gamma(\IP(V_\eta),\O_{S^d(\X/B)}(\H_{s,t})|_{\IP(V_\eta)})=
\binom{s+t+r}{r}
$$
for $s,t>>0$. On the other hand, also
$$
P_{\IP(\mathfrak g)}(s,t):=
\dim\Gamma(\IP(\mathfrak g),\O_{S^d(\X/B)}(\H_{s,t})|_{\IP(\mathfrak g)})=
\binom{s+t+r}{r}
$$
for $s,t>>0$ by Theorem \ref{main}. Thus, we need only show that the 
closure of $\IP(V_\eta)$ contains $\IP(\mathfrak g)$. 

Consider now on the product $\X\times_B\IP(\wt V_i)$ the 
composition
$$
\O_{\IP(\wt V_i)}(-1) \lra \wt V_i \lra \L(-iZ)
$$
where the first map is the tautological map of $\IP(\wt V_i)$ and the 
second is the evaluation map, all sheaves 
and maps
being viewed on the product 
under the appropriate pullbacks. Let $F$ denote the degeneracy scheme 
of this composition. Then $F$ is a relative Cartier divisor of degree $d$ 
of $\X\times_B\IP(\wt V_i)/\IP(\wt V_i)$ away from 
$\IP(V_i^{Y,0})\cup\IP(V_i^{Z,0})$. Thus we obtain a map 
$$
\IP(\wt V_i)-(\IP(V_i^{Y,0})\cup\IP(V_i^{Z,0}))\lra S^d(\X/B)
$$
whose image contains $\IP(V_\eta)$ and all points of $S^d(X)$ of the 
form $\div(s|_Y)+\div(s|_Z)$ for $s\in V_i-(V_i^{Y,0}\cup V_i^{Z,0})$. 
Since $\IP(\wt V_i)$ is flat over $B$, it follows that the closure 
of $\IP(V_\eta)$ in $S^d(X)$ contains all points of the above form. As 
we let $i$ vary, we get that the closure of $\IP(V_\eta)$ contains 
$\IP(\mathfrak g)$.
\end{proof}

\section{Limitations on limit linear series}

In view of Theorem \ref{main}, if a given fiber $A_d^{-1}(\L)$ of the Abel 
map has dimension strictly less than $r$, it is clear that there cannot 
exist any exact limit linear series on $X$ of dimension $r$ having $\L$ as
its underlying line bundle. However, a substantially stronger assertion
holds. Namely, we have:

\begin{prop}\label{prop:no-grds} If $A_d^{-1}(\L)$ has any irreducible 
component of 
dimension
strictly smaller than $r$, then there is no limit linear series on $X$
of dimension $r$ and underlying line bundle $\L$.
\end{prop}

For the proof, it will be convenient to work with Eisenbud--Harris limit
linear series, so we briefly recall the definition. 

\begin{defn} An \emph{Eisenbud--Harris limit linear series} on $X$ of
dimension $r$ and degree $d$ consists of a pair $((\L^Y,V^Y),(\L^Z,V^Z))$
of linear series on $Y$ and $Z$, each of dimension $r$ and degree $d$,
such that if $a_0^Y,\dots,a_r^Y$ and $a_0^Z,\dots,z_r^Z$ denote the
vanishing sequences of $(\L^Y,V^Y)$ and $(\L^Z,V^Z)$ respectively at $P$,
we have
$$a_i^Y+a_{r-i}^Z \geq d$$
for each $i=0,\dots,r$. 
\end{defn}

Given a limit linear series $(\L,V_0,\dots,V_d)$ in our sense, we clearly
obtain a pair of linear series by setting $\L^Y=\L^0|_Y$,
$\L^Z=\L^d|_Z$, $V^Y=V_0|_Y$, $V^Z=V_d|_Z$. According to Proposition 6.6
of \cite{os8}, the resulting pair is in fact an Eisenbud--Harris limit
linear series, and every Eisenbud--Harris limit linear series arises in
this way. 
As a consequence,
the statement of Proposition \ref{prop:no-grds} is equally valid in either 
context. Note also that 
the converse of the statement does not hold: 
even if the fiber of the Abel maps has every
component of dimension at least $r$, there may be no limit linear series
of dimension $r$; see Example \ref{ex:no-grds} below.

We will also find the following terminology convenient:

\begin{defn} Suppose $\L$ is a line bundle on $X$ of degree $d$ on $Y$ 
and degree $0$ on 
$Z$. 
Then the \emph{vanishing sequence} of $\L$ on $Y$
at $P$ is the vanishing sequence at $P$ of the complete linear series for
$\L|_Y$. Similarly, the \emph{vanishing sequence} of $\L$ on $Z$ at $P$
is the vanishing sequence at $P$ of the complete linear series for
$(\L|_Z)(dP)$.
\end{defn}

\begin{lemma}\label{lem:fiber-comps} Let $\L$ be a line bundle 
on $X$ of degree $d$ on $Y$ and degree $0$ on $Z$, 
with vanishing sequences $a_0^Y,\dots,a_p^Y$ and $a_0^Z,\dots,a_q^Z$ on
$Y$ and $Z$ respectively at $P$. For $\ell$ between $0$ and $d$, the 
subset $\IP(\Gamma^\ell_Y) \times \IP(\Gamma^{d-\ell}_Z)$ inside
$A_d^{-1}(\L)$ constitutes an irreducible component of the fiber if and 
only if there exist $\ell' \leq \ell, \ell'' \geq \ell,$ and $i,j$ such 
that 
$a_i^Y=d-\ell'$, $a_j^Z=\ell''$, $a_{i-1}^Y < d-\ell''$, and $a_{j-1}^Z<\ell'$.
In this case, the dimension of 
$\IP(\Gamma^\ell_Y) \times \IP(\Gamma^{d-\ell}_Z)$ is equal to
$p+q-i-j$. 
\end{lemma}

\begin{proof} 
Because of \eqref{eq:abel-fiber}, since each 
subset $\IP(\Gamma^\ell_Y) \times \IP(\Gamma^{d-\ell}_Z)$ is irreducible, 
one such subset is an irreducible component 
of the fiber if and only if it is not strictly contained in any other. There
are only two ways such a strict containment could occur: 
as $\ell$ decreases, we could have 
$\IP(\Gamma^{\ell}_Y)$ remaining unchanged until after 
$\IP(\Gamma^{d-\ell}_Z)$ has strictly increased, or as $\ell$ increases, 
we could have $\IP(\Gamma^{d-\ell}_Z)$ remaining unchanged until
after $\IP(\Gamma^{\ell}_Y)$ has strictly increased.
These two conditions
are precisely what is ruled out by the conditions in the statement of the
lemma.
\end{proof}

\begin{proof}[Proof of Proposition \ref{prop:no-grds}] According to 
Lemma \ref{lem:fiber-comps}, any irreducible component of the fiber
is of dimension $p+q-i-j$, where $i,j$ satisfy $a_{i-1}^Y+a_j^Z<d$. 
We prove that if there is an Eisenbud--Harris limit
linear series with underlying line bundle $\L$ and dimension $r$, then
we must have $p+q-i-j \geq r$. Given
$\L$, such a limit series is determined by the spaces $V^Y,V^Z$, which
must have their vanishing sequences at $P$ being subsequences of length
$r+1$ of the $a_s^Y$ and $a_s^Z$; suppose these are given by indices
$i_0<i_1 < \dots < i_r$ and $j_0 < j_1 <\dots < j_r$. Then by the
definition of an Eisenbud--Harris limit series,
we must have $a_{i_s}^Y+a_{j_{r-s}}^Z \geq d$ for $s=0,\dots,r$.  
Choose $s$ maximal with $i_s \leq i-1$. Then $i_{s+1},\dots,i_r \geq i$,
so in order to have enough remaining choices of indices, we conclude that 
$r-s \leq p+1-i$. On the other hand, since
$a_{i_s}^Y+a_{j_{r-s}}^Z \geq d$, we have
$$a_{j_{r-s}}^Z \geq d-a_{i_s}^Y \geq d-a_{i-1}^Y>a_j^Z,$$ 
so as before we conclude 
$s+1=r+1-(r-s) < q+1-j$. Putting 
the two
inequalities together we conclude that $r \leq p+q-i-j$, as desired.
\end{proof}

\section{Examples and further discussion}

In general, the map from exact limit linear series to subschemes
of fibers of Abel maps need not be injective. The reason is essentially
that if $V_i$ is obtained by gluing together sections of 
$\IP(\Gamma^{d-i}_Y)$ and $\IP(\Gamma^i_Z)$ 
which vanish at the node $P$, but we do not
have $V_i=V_i^{Y,0} \oplus V_i^{Z,0}$, then scaling sections on $Z$ while
holding them fixed on $Y$ will yield different choices for the subspaces
$V_i$, but will not change the associated subsets of $S^d(X)$. For instance,
if $r=0$ an exact limit linear series is uniquely determined by a single 
choice of section of some $\L^i$ which does not vanish on either $Y$ or $Z$. 
If this section vanishes at $P$, we may obtain different choices of $V_i$
by scaling the section on $Z$ while holding its value on $Y$ fixed. In this
case, we see that the map from the moduli space of 
exact
limit linear series to
that of 
subschemes of fibers of the Abel map actually factors through the
space of 
Eisenbud--Harris limit linear series. However, this is not the case in
general.

\begin{ex} Suppose $X$ has genus $0$, and 
set $d=2$, $r=1$. In this
case we have a $2$-dimensional family of exact limit series which all
correspond to the same Eisenbud--Harris limit series. Explicitly, if
we choose coordinates $y,z$ on $Y$ and $Z$ such that $P$ is $y=0$ and $z=0$,
we may represent sections by pairs of polynomials in $y$ and $z$ of degree
at most $2$. Then $(y,y^2), (z,z^2)$ gives a (crude) Eisenbud--Harris limit
linear series. The corresponding choices of limit linear series in our sense
are (as always) uniquely determined for $V_0,V_2$, but $V_1$ may be
any $2$-dimensional subspace of the vector space of pairs 
$(a_1 y+a_2 y^2,a_1 z+a_3 z^2)$. 

This vector space is $3$-dimensional,
so we have a $2$-dimensional projective space of choices for $V_1$. The
exact limit series are an open subset (specifically, those spaces 
which remain $2$-dimensional after restriction to either $Y$ or $Z$), and
we see that depending on our choice of $2$-dimensional subspace, we will
obtain different pencils of corresponding divisors, and hence different
$1$-dimensional subschemes of the fiber of the Abel map.
\end{ex}

Thus, we cannot define our map using Eisenbud--Harris limit series instead
of our limit series, and in particular this is not a viable approach to
resolving the lack of injectivity. However,
it appears that the failure of injectivity may be resolved if instead of
associating subschemes of $S^d(X)$ to a given limit series, we associate
subschemes of $\Hilb^d(X)$. 

One may also wonder to what extent the map
we have constructed is surjective. There is more than one way to interpret
this question, and we consider more specifically whether each fiber of the
Abel map is the union of the closed subschemes we construct. For this
question, it is necessary to specify further which $r$ we consider. If we
take $r=0$, we have very few constraints on constructing limit linear
series, but since every exact limit series yields a point of the fiber
of the Abel map, and since we do not associate subschemes of the fiber
to a non-exact limit linear series, we should only expect to get an
open subset of the fiber in this way, and this is indeed what happens.

\begin{ex} Suppose $Y$ and $Z$ are both elliptic curves. Consider
the case $d=1$, and $\L$ the line bundle of degree $1$ such that the
degree-$1$ restrictions to $Y$ and $Z$ are both isomorphic to $\O(P)$.
Then each 
restriction has a unique nonzero 
section up to scalar, which
vanishes to order $1$ at $P$. The fiber $A^{-1}_1(\L)$ is the single
point $P$, which we may view either as 
$\IP(\Gamma_Y^0) \times \IP(\Gamma_Z^1)$ or as 
$\IP(\Gamma_Y^1) \times \IP(\Gamma_Z^0)$. 
However, while there is a
non-exact limit series of dimension $0$ with underlying line bundle 
$\L$, there is no exact one, so we are unable to associate any subschemes
to this fiber.
\end{ex}

The lack of surjectivity exhibited in the previous example could in
principle be addressed by extending our construction to non-exact limit
linear series. However, we see that for $r>0$, there may not be any
limit linear series even when the fiber of the Abel map is pure of
dimension $r$.

\begin{ex}\label{ex:no-grds} With $Y$ and $Z$ still both elliptic curves, 
suppose $d=2$,
and the degree-$2$ restrictions of $\L$ to $Y$ and to $Z$ are 
$\O_Y(2P)$ and $\O_Z(P+Q)$ for some $Q \neq P$, respectively.
Then $\L$ has vanishing sequences at $P$ given by $a_0^Y=0,a_1^Y=2$ and 
$a_0^Z=0, a_1^Z=1$. We thus see that the corresponding fiber of $A_2$
is equal to 
$\IP(\Gamma_Y^0) \times \IP(\Gamma_Z^2)\cong \IP^1$.
However, we see immediately that there is no Eisenbud--Harris limit
linear series of dimension $1$ with underlying line bundle $\L$,
and thus no limit linear series in our sense, either.
\end{ex}

Thus, to address the surjectivity question one also has to determine what 
range of $r$ is appropriate to consider for a given fiber of the Abel map.

Finally, we give an example for which the fiber of the Abel map is not
equidimensional. 

\begin{ex} With $Y$ and $Z$ still both elliptic curves, suppose $d=3$,
and the degree-$3$ restrictions of $\L$ to $Y$ and $Z$ are
$\O_Y(3P)$ and $\O_Z(2P+Q)$ for some $Q \neq P$, respectively. Then $\L$
has vanishing sequences at $P$ given by $a_0^Y=0,a_1^Y=1,a_2^Y=3$, and
$a_0^Z=0, a_1^Z=1, a_2^Z=2$. In this case, we see that $A_3^{-1}(\L)$
is the union of 
$\IP(\Gamma_Y^0)\times\IP(\Gamma^3_Z)$ and 
$\IP(\Gamma_Y^2) \times \IP(\Gamma^1_Z)$, which have dimension $2$ and $1$ 
respectively, and intersect only at a point.
\end{ex}

\bibliographystyle{hamsplain}
\bibliography{hgen}

\end{document}